\documentclass[a4paper,12pt]{article}
\usepackage{amsthm,amsmath,amssymb}
\newtheorem{teor}{Theorem}[section]
\newtheorem{defin}[teor]{Definition}
\newtheorem{lemm}[teor]{Lemma}
\newtheorem{osse}[teor]{Remark}
\newtheorem{prop}[teor]{Proposition}
\newtheorem{defi}[teor]{Definition}
\newtheorem{coro}[teor]{Corollary}
\newtheorem{prob}[teor]{Problem}

\newcommand{\bele}{\begin{lemm}\begin{sl}}
\newcommand{\enle}{\end{sl}\end{lemm}}
\newcommand{\bedef}{\begin{defi}\begin{sl}}
\newcommand{\eddef}{\end{sl}\end{defi}}
\newcommand{\bete}{\begin{teor}\begin{sl}}
\newcommand{\ente}{\end{sl}\end{teor}}
\newcommand{\beos}{\begin{osse}\begin{rm}}
\newcommand{\eddos}{\end{rm}\end{osse}}
\newcommand{\bepr}{\begin{prop}\begin{sl}}
\newcommand{\empr}{\end{sl}\end{prop}}
\newcommand{\bepro}{\begin{prob}\begin{rm}}
\newcommand{\empro}{\end{rm}\end{prob}}
\newcommand{\bede}{\begin{defin}\begin{sl}}
\newcommand{\edde}{\end{sl}\end{defin}}
\newcommand{\beco}{\begin{coro}\begin{sl}}
\newcommand{\enco}{\end{sl}\end{coro}}

\newcommand{\quand}{\quad\text{and}\quad}

\newcommand{\quext}{\quad\text}

\newcommand{\de}{\partial}

\newcommand{\SSS}{\mathbb{S}}
\newcommand{\RR}{\mathbb{R}}
\newcommand{\II}{\mathbb{I}}
\newcommand{\EE}{\mathbb{E}}
\newcommand{\NN}{\mathbb{N}}

\newcommand{\beeq}[1]{\begin{equation}\label{#1}}
\newcommand{\eddeq}{\end{equation}}

\newcommand{\beeqa}[1]{\begin{eqnarray}\label{#1}}
\newcommand{\eddeqa}{\end{eqnarray}}

\newcommand{\beal}[1]{\begin{align}\label{#1}}
\newcommand{\eddal}{\end{align}}

\newcommand{\bespl}[1]{\begin{split}\label{#1}}
\newcommand{\edspl}{\end{split}}

\newcommand{\bega}[1]{\begin{gather}\label{#1}}
\newcommand{\edga}{\end{gather}}

\newcommand{\beeqax}{\begin{eqnarray*}}
\newcommand{\eddeqax}{\end{eqnarray*}}

\def\qed{\ifmmode 
  \else \leavevmode\unskip\penalty9999 \hbox{}\nobreak\hfill
  \fi
  \quad\hbox{\hskip.5em\vrule width.4em height.6em depth.05em\hskip.1em}}
\def\endproofsym{\qed}
\renewenvironment{proof}[1][Proof]{\trivlist\item[\hskip\labelsep{\hskip0pt
    {\normalfont\scshape#1.}\hskip .321429\parindent}]\ignorespaces}
{\endproofsym\endtrivlist}
\def\endnobox{\def\endproofsym{}\end{proof}\def\endproofsym{\qed}}

\newcommand{\no}{\nonumber}

\newcommand{\beeqao}{\begin{eqnarray}\no}
\newcommand{\bealo}{\begin{align}\no}
\newcommand{\besplo}{\begin{split}\no}
\newcommand{\begao}{\begin{gather}\no}

\newcommand{\duav}[1]{\langle{#1}\rangle}

\newcommand{\perogni}{\forall\,}

\newcommand{\io}{\int_\Omega}

\newcommand{\OO}{_{\Omega}}

\newcommand{\psihat}{\widehat{\psi}}

\newcommand{\bu}{\boldsymbol{u}}

\newcommand{\di}{\boldsymbol{d}}
\newcommand{\ddi}{\overline{\di}}
\newcommand{\dii}{\di_{\infty}}

\newcommand{\pp}{\boldsymbol{p}}

\newcommand{\bz}{\boldsymbol{z}}

\newcommand{\bV}{\boldsymbol{V}}
\newcommand{\bVz}{\boldsymbol{V}_{\mbox{}\!0}}
\newcommand{\bH}{\boldsymbol{H}}

\newcommand{\bh}{\boldsymbol{h}}

\newcommand{\bd}{\boldsymbol{d}}

\newcommand{\bn}{\boldsymbol{n}}

\newcommand{\bv}{\boldsymbol{v}}
\newcommand{\dn}{\partial_{\bn}}

\newcommand{\bbf}{\boldsymbol{f}}

\newcommand{\bphi}{\boldsymbol{\phi}}

\newcommand{\rhs}{right-hand side}

\DeclareMathOperator{\dive}{div}

\DeclareMathOperator{\deriv}{d}

\let\TeXchi\chi
\def\chi{{\setbox0 \hbox{\mathsurround0pt
$\TeXchi$}\hbox{\raise\dp0 \copy0 }}}

\newcommand{\calD}{{\mathcal D}}

\newcommand{\calE}{{\mathcal E}}

\newcommand{\calM}{{\mathcal M}}

\newcommand{\calB}{{\mathcal B}}

\newcommand{\calV}{{\mathcal V}}

\newcommand{\dit}{\deriv\!t}
\newcommand{\dis}{\deriv\!s}

\newcommand{\ddt}{\frac{\deriv\!{}}{\dit}}

\newcommand{\ii}{_\infty}

\newcommand{\calL}{{\cal L}}

\begin{document}

\title{On the long-time behavior of some mathematical models
for nematic liquid crystals}

\author{Hana Petzeltov\'a\\
{\sl Institute of Mathematics of the Czech Academy of Sciences}\\
{\sl \v Zitn\' a 25, 115 67 Praha 1}\\
{\sl Czech Republic}\\
{\rm E-mail:~~\tt petzelt@math.cas.cz}
\and
Elisabetta Rocca\thanks{The work of E.R.~was supported
by the FP7-IDEAS-ERC-StG Grant \#256872 (EntroPhase)}\\
{\sl Dipartimento di Matematica, Universit\`a di Milano}\\
{\sl Via Saldini, 50, I-20133 Milano, Italy}\\
{\rm E-mail:~~\tt elisabetta.rocca@unimi.it}
\and
Giulio Schimperna\thanks{The work of G.S.~was partially supported by
Grant as a part of the general research programme of the Academy of
Sciences of the Czech Republic, Institutional Research Plan
AV0Z10190503}\\
{\sl Dipartimento di Matematica, Universit\`a di Pavia}\\
{\sl Via Ferrata, 1, I-27100 Pavia, Italy}\\
{\rm E-mail:~~\tt giusch04@unipv.it} }
\date{}
\maketitle

\begin{abstract}
 A model describing the evolution of a liquid crystal
 substance in the nematic phase is investigated in terms of two
 basic state variables: the {\it velocity field} $\bu$
 and the {\it director field} $\di$, representing the
 preferred orientation of molecules in a neighborhood of any point
 in a reference domain. After recalling a known existence result,
 we investigate the long-time behavior of weak solutions.
 In particular, we show that any solution trajectory
 admits a non-empty $\omega$-limit set containing only
 stationary solutions. Moreover, we give a number of
 sufficient conditions in order that the $\omega$-limit
 set contains a single point. Our approach improves
 and generalizes existing results
 on the same problem.
\end{abstract}

\medskip


\medskip

\noindent
{\bf AMS (MOS) subject clas\-si\-fi\-ca\-tion:}~~35B40, 35K45, 76A15.

\maketitle


\section{Introduction}
\label{sec:intro}

In this paper we analyze the long-time behavior of weak solutions to
the system
\begin{align}\label{eq:uin}
  & \bu_t + \dive (\bu \otimes \bu)
   - \nu \Delta \bu = \dive \big( - p \II
   - L (\nabla \di \odot \nabla \di)
   - \delta (L \Delta \di - \bbf(\di)) \otimes \di \big),\\
 \label{incomprin}
  & \dive \bu = 0,\\
 \label{eq:din}
  & \di_t + \bu \cdot \nabla \di
   - \delta \di \cdot \nabla \bu
   - L \Delta \di + \bbf(\di) = 0,
\end{align}
describing the evolutionary behavior of {\sl nematic}\/ liquid
crystal flows (we refer to the monographs \cite{Ch,dG} for
a detailed presentation of the physical foundations
of continuum theories of liquid crystals). Actually,
system \eqref{eq:uin}-\eqref{eq:din}
can be seen as a simplification of the original
Ericksen-Leslie model \cite{eric,le},
that still keeps a good level of compliance
with experimental results. The model couples the
Navier-Stokes equation~\eqref{eq:uin}
for the macroscopic velocity $\bu$
($p$ denoting as usual the pressure),
with the incompressibility condition
\eqref{incomprin} and with the equation \eqref{eq:din}
ruling the behavior of the local orientation
vector $\di$ of the liquid crystal. Here, the function
$\bbf$ represents the gradient w.r.t.~$\di$
of the configuration energy $F$ of the crystal.
We choose $F$ to be a double well potential
having minima for $|\di|=1$ and growing at infinity
at most as a fourth order polynomial. This provides
a standard relaxation of the physical constraint
$|\bd|=1$, which is very difficult to treat mathematically.

In this paper, the system is
complemented with the homogeneous Dirichlet boundary
condition for $\bu$, the no-flux condition for $\di$, and
with initial conditions. It is settled in a smooth
bounded domain $\Omega\subset \RR^d$ for $d=2$ or $d=3$.
No restriction is assumed on the viscosity
coefficient $\nu$.

Regarding the parameter $\delta$, we
will take $\delta\ge 0$, with the case $\delta>0$ denoting
the presence of a {\sl stretching effect}\/
on the molecules of the crystal. Some of our results,
however, hold only for $\delta=0$.
Actually, the situation $\delta>0$
is more difficult to be
treated mathematically since
the term $\delta\di \cdot \nabla \bu$ prevents
from using maximum principle arguments in
\eqref{eq:din}. For this reason, even if the
initial datum $\di_0$ satisfies the (relaxed)
physical constraint $|\di_0|\le 1$ almost
everywhere,  the same may not be true
for $\di(t)$, for positive times, if $\delta>0$.

A mathematical analysis of system \eqref{eq:uin}-\eqref{eq:din}
has been first addressed in the papers \cite{LinLiusimply}
and \cite{LinLiu} (in this second work, an even more general model
is taken into account). There, the authors consider
the case $\delta=0$ and prove existence of a unique classical
solution for $d=2$, and also in dimension $d=3$ under the additional
assumption that the viscosity $\nu$ is sufficiently large.
These results have been extended to
the case $\delta>0$ in the paper \cite{SunLiu}.
Finally, the restriction on the viscosity has been
recently dropped in \cite{CR}, where {\sl weak}\/
solutions are considered and a global
existence result for the 3D system
\eqref{eq:uin}-\eqref{eq:din} is proved
in that regularity frame. Of course,
uniqueness is not known to hold in that regularity setting.
A similar result is essentially contained also in the
recent paper \cite{CGR}, where analogous estimates
are derived but no formal statement of an existence
result is provided.

The Dirichlet boundary condition for $\bu$ and
either a nonhomogeneous Dirichlet or the no-flux
boundary condition for $\di$ are treated there. Moreover, let
us quote the recent paper \cite{FRS2}, where
these results have been extended to a more general
system \eqref{eq:uin}-\eqref{eq:din}, where also
temperature effects are taken into account.
We note, however, that the results of \cite{FRS2}
require different boundary conditions for $\bu$ (namely,
the so-called {\sl complete slip}\/ conditions).

The long-time behavior of system
\eqref{eq:uin}-\eqref{eq:din}
has been analyzed in the recent work \cite{WXL},
still considering the case $d=2$ or the case
$d=3$ with the large viscosity $\nu$, and periodic boundary conditions. More precisely, in \cite{WXL}
the authors show existence of a nonempty $\omega$-limit set
for any strong bounded solution emanating from smooth initial data. Moreover, by using the
Simon-\L ojasiewicz inequality, they prove that,
for the nonlinearity $\bbf =(|\di|^2-1)\di$, this
$\omega$-limit set contains only one point.

Stability and asymptotic stability properties of this model
(actually, with even more complete
stretching terms) have also been studied in \cite{CGR}, where the long-time behavior of
solutions is analyzed in the case of periodic boundary conditions.
More precisely, the authors prove, by means of formal estimates,
that weak solutions become eventually smoother for large times,
which suffices to have existence of non-empty $\omega$-limit sets.

Finally, in the recent contribution \cite{GrasWu},
the existence of a smooth global attractor of finite fractal dimension is obtained in
two dimensions of space.

Our aim in this paper is to extend the results
of \cite{CGR,WXL} in the following directions:\\[2mm]
{\sl (i)}~~we address the case $d=3$ without
the large viscosity assumption considering weak solutions;\\[1mm]
{\sl (ii)}~~we consider more general $C^1$ functions $\bbf$; \\[2mm]
{\sl(iii)}~~we use different boundary conditions and weaker initial data; \\[2mm]
{\sl(iv)}~~we discuss convergence, as $t$ tends to $\infty$, of strong solutions in some particular situations.\\[2mm]
To get {\sl (i)}, we prove convergence of weak solutions, and,
in some situations, we get strong convergence using the fact
that weak solutions to the system become eventually
smoother for times $t$ larger than some $T$.
This property is well-known for the (uncoupled)
three-dimensional N-S system, and we find conditions under which
it holds also for the coupled system
\eqref{eq:uin}-\eqref{eq:din}. Note that this result is still true for periodic boundary conditions,
and so it improves the study done  in \cite{WXL}. In turn, this property
(cf.~\eqref{couv}-\eqref{regoult:d} below)
enables us to obtain properties sufficient
to characterize the $\omega$-limit set. Assuming that $\bbf$ is analytic, we
apply the generalized \L ojasiewicz theorem to get convergence of the variable $\di$.

To address question {\sl (ii)}, in particular, to remove the analyticity condition,
we make the basic
observation that the set of global minimizers of the
configuration energy of the crystal coincides
with the set of constant unit vectors of $\RR^d$.
Then, it is easy to prove
that any global minimizer $\ddi$ satisfies the so-called
{\sl normal hyperbolicity}\ condition.  Based on this fact, we can prove
that, if the $\omega$-limit set contains
a global minimizer, then it coincides with
it (i.e., it does not contain any other point).
We can also give two precise conditions ensuring
the fact that the $\omega$-limit set contains global minimizers, which, unfortunately, require
$\delta=0$. Namely, this happens when
either the diffusion coefficient $L$
is large enough, or when
the initial energy is very small compared with $L$ (in particular,
the initial datum $\di_0$ is already close enough to the set of global minimizers
in a suitable norm).

The paper is organized as follows. In the next
section, we present our assumptions, state the main
results, and, for the reader's convenience, we
briefly sketch the basic estimates at the core
of the existence proof. The proofs of the new
results on the long-time behavior are given in
Section~\ref{sec:proofs}.

\medskip

{\bf Acknowledgment.}~~The authors are grateful to Professor Eduard Feireisl
 for illuminating discussions regarding the results proved in this paper.


\section{Main results}
\label{sec:results}

We let $\Omega$ be a smooth, bounded, and connected
domain in $\RR^d$, $d\in\{2,3\}$, with the boundary $\Gamma$.
For simplicity, we also assume $|\Omega|=1$.
We set $H:=L^2(\Omega)$, $\bH:=L^2(\Omega)^d$, and
denote by $(\cdot,\cdot)$ the scalar product
both in $H$ and in $\bH$ and by $\|\cdot\|$
the related norms.
Next, we set $V:=H^1(\Omega)$, $\bV:=H^1(\Omega)^d$
and $\bVz:=H^1_0(\Omega)^d$.
The duality between $V'$ and $V$,
as well as those between $\bV'$ and $\bV$
and between $\bVz'$ and $\bVz$,
will be indicated by $\langle\cdot,\cdot\rangle$.
Identifying $H$ with $H'$ through the scalar product
of $H$, it is then well known that
$V\subset H\subset V'$ with continuous
and dense inclusions. In other words,
$(V,H,V')$ constitutes a Hilbert triplet
(see, e.g., \cite{Li}). Correspondingly, we
also have the vectorial analogues $(\bV,\bH,\bV')$
and $(\bVz,\bH,\bVz')$. The symbol $\|\cdot\|_{X}$
will indicate the norm
in the generic (real) Banach space $X$ and
$\langle\cdot,\cdot\rangle_X$ will
stand for the duality between $X'$ and $X$.

\smallskip

We consider $\bbf$ in the form
\begin{equation}\label{defibbf}
   \bbf(\di)=\big(\psi(|\di|^2)-1\big)\di=\frac12\partial_{\di}\big(\psihat(|\di|^2)-|\di|^2\big),
\end{equation}
where
\begin{equation}\label{hp:psi}
  \psi\in C^1\big([0,+\infty);[0,+\infty)\big),
   \quext{with }\,\psi(0)=0,~~\psi(1)=1~~
   \text{and }\,\psi'(1)>0,
\end{equation}
is an increasing function, and the convex function $\psihat$ is defined by
\begin{equation}\label{hp:psihat}
  \psihat'=\psi,\quad \psihat(1)=1.
\end{equation}
We also assume that there exists a constant $c_\psi>0$ such that
\begin{equation}\label{hp:psibound}
  \psi'(r) \le c_\psi
   \quext{for all~}\ r\in [0,+\infty).
\end{equation}
Given $L>0$, we define the
{\sl configuration energy}\ of the
liquid crystal flow as
\begin{equation}\label{defiE}
  \calE(\di):= \frac12 \io
   \big(L |\nabla \di|^2
    + \psihat(|\di|^2)
    - |\di|^2 \big).
\end{equation}
The {\sl total energy}\ is then given by
adding to $\calE$ the ``macroscopic''
kinetic energy; namely, we set
\begin{equation}\label{defiEE}
  \EE(\bu,\di):= \frac12 \| \bu \|^2 + \calE(\di)
   = \frac12 \io
   \big(| \bu |^2
    + L |\nabla \di|^2
    + \psihat(|\di|^2)
    - |\di|^2 \big).
\end{equation}
Let us notice that, thanks to the above
assumptions \eqref{hp:psi}-\eqref{hp:psibound},
$\calE(\di)=0$ if and only if $\di$ is a (constant)
unit vector (cf. Lemma~\ref{lemma:gm} below for a simple proof).

\smallskip

We will address the following system of PDE's:
\begin{align}\label{eq:u}
  & \bu_t + \dive (\bu \otimes \bu)
   - \nu \Delta \bu = \dive \SSS,\\
 \label{eq:S}
  & \SSS = - p \II
   - L (\nabla \di \odot \nabla \di)
   - \delta (L \Delta \di - \bbf(\di)) \otimes \di,\\
 \label{incompr}
  & \dive \bu = 0,\\
 \label{eq:d}
  & \di_t + \bu \cdot \nabla \di
   - \delta \di \cdot \nabla \bu
   - L \Delta \di + \bbf(\di) = 0,
\end{align}
where  the coefficients $\nu,L,\delta$ satisfy
$\nu,L>0$ and $\delta\ge0$.
Notice that, by \eqref{hp:psibound}, $\bbf(\di)$ grows
at infinity at most as the third power of $|\di|$.

The system, supplemented  with the boundary and initial conditions
\begin{align}\label{bc:u}
  & \bu = 0 \quext{a.e.~on }\,
   (0,T)\times \Gamma,\\
 \label{bc:d}
  & \dn \bd = 0 \quext{a.e.~on }\,
   (0,T)\times \Gamma,
\end{align}
\begin{equation}\label{init}
  \bu|_{t=0}=\bu_0, \quad
   \di|_{t=0}=\di_0,
    \quext{a.e.~in }\,\Omega,
\end{equation}
will be called Problem (P).

We introduce a precise definition
of weak solutions:
\bede\label{def:weaksol}
  {\rm  A weak solution} to\/ {\rm Problem (P)} is a couple $(\bu,\di)$
  such that
 \begin{align}\label{reg:u}
   & \bu \in L^\infty(0,T;\bH) \cap L^2(0,T;\bVz),\\
  \label{reg:d}
   & \di\in H^1(0,T;L^{3/2}(\Omega)^d)\cap L^\infty(0,T;\bV)
    \cap L^2(0,T;H^2(\Omega)^d),
 \end{align}
 for all $T>0$, $\bu,\di$ satisfy initial and boundary conditions \eqref{init}, \eqref{bc:d},
 the equations \eqref{eq:S}-\eqref{eq:d} are satisfied
 for a.e.~$t\in(0,T)$, and
 \begin{equation}\label{eq:u-w}
   \duav{\bu_t,\bphi} - \io (\bu \otimes \bu) : \nabla\bphi
    + \nu \io \nabla \bu : \nabla \bphi = - \io \SSS : \nabla \bphi,
 \end{equation}
 holds for any test function $\bphi \in \boldsymbol{W}^{1,3}_{0,\dive}(\Omega)$
 (i.e., the subspace of $W^{1,3}_0(\Omega)^d$ consisting
 of divergence-free functions).
\edde
\beos\label{rem:weaksol}
 The regularity of the test function $\phi$ can be justified thanks to \eqref{reg:u}, \eqref{reg:d} and \eqref{eq:d}. We have
 in any case (also if $\delta>0$)
 \begin{equation}\label{on:reg}
  \bu\otimes \bu,~~
   \nabla \di \odot \nabla \di,~~
    (L \Delta \di - \bbf(\di)) \otimes \di
   \in L^2(0,T;L^{3/2}(\Omega))^{d\times d},
 \end{equation}
 whence their (distributional) divergence belongs to
 the space $L^2(0,T;W^{-1,3/2}(\Omega))^{d}$.
 Note also that the boundary condition~\eqref{bc:u}
 is in fact ``embedded'' into the weak formulation \eqref{eq:u-w}.
\eddos
\noindent%
It is known that Problem (P) admits
at least one weak solution $(\bu,\di)$.
This has been proved in \cite{LinLiusimply}
for the case $\delta=0$ and in \cite{CR}
for the case $\delta=1$ (cf. also \cite{CGR} for the formal computations).
Namely, we have
\bete\label{teo:esi}
 Let\/ \eqref{defibbf}-\eqref{hp:psibound} hold and let
 \begin{gather}\label{hp:init1}
   \bu_0 \in \bH, \quad
    \dive \bu_0 = 0,\\
  \label{hp:init2}
   \bd_0 \in \bV.
 \end{gather}
 Then, {\rm Problem (P)} possesses a global in time weak solution
 $(\bu,\di)$, satisfying,
  for a.a.~$t>0$, the\/
 {\rm energy inequality}
 \begin{equation}\label{energy}
   \ddt \EE(\bu,\di)
    + \big\| - L \Delta \di + \bbf(\di) \big\|^2
    + \nu \| \nabla \bu \|^2 \le 0.
 \end{equation}
 \ente
\noindent%
We point out that assumptions
\eqref{hp:init1}-\eqref{hp:init2} are equivalent to asking
that the initial energy
$\EE_0:=\EE(\bu_0,\di_0)$ is finite.

\beos\label{which-sol}
 The proof of the above theorem relies on a
 rather tricky approximation scheme and on refined
 compactness methods to pass to the limit.
 It is then worth pointing out that, due to
 nonuniqueness, our subsequent results
 on the long-time behavior
 hold only for those solutions satisfying the energy inequality, in particular for
 the limit points of the approximate scheme,
 and not necessarily for all
 solutions in the regularity frame
 \eqref{reg:u}-\eqref{reg:d}.
 Actually, there may
 exist ``spurious'' weak solutions not satisfying
 the energy inequality \eqref{energy} which is crucial for investigating
  the long-time behavior. As a convention,
 in the sequel we shall restrict the terminology
 ``weak solutions'' to those solutions which satisfy \eqref{energy}. Spurious
 solutions are thus excluded.
\eddos

As noted above, in the case $\delta=0$
a maximum principle holds for
the $\di$-component of any weak solution.
For the reader's convenience, we recall the statement
and the (simple) proof.
\bete\label{teo:max}
 Let the assumptions of\/ {\rm Theorem~\ref{teo:esi}} hold
 and let
 \[
 \delta=0, \mbox{ and}\quad  |\di_0(x)|\le 1 \mbox{ for a.a. }x\in\Omega.
 \]
  Then
 any weak solution $(\bu,\di)$ to\/ {\rm Problem (P)} satisfies
 \begin{equation}\label{maxprinc}
   |\di(t,x)|\le 1  \quext{for a.a.~}\,
    (t,x)\in (0,\infty)\times \Omega.
 \end{equation}
\ente
\begin{proof}
 Testing equation \eqref{eq:d} by $\di$ one obtains
 \begin{equation}\label{d:scalar}
   \frac12 \ddt |\di|^2 + \bu \cdot \nabla | \di |^2
    - \frac{L}2 \Delta |\di|^2 + L | \nabla \di |^2
    + \big( \psi(|\di|^2) - 1\big)|\di|^2 = 0.
 \end{equation}
 Then, we notice that, by \eqref{hp:psi},
 \begin{equation}\label{f:maxprinc}
   (\psi(r)-1)r \ge 0 \quad \perogni r\ge 1.
 \end{equation}
 Thus, \eqref{d:scalar} represents a parabolic
 equation for $|\di|^2$ (which still satisfies the
 no-flux b.c.). It is clear that the maximum principle
 applies, yielding \eqref{maxprinc}.
\end{proof}
\noindent%
Unfortunately, \eqref{maxprinc} is not known (and not expected)
to hold in the case $\delta>0$.

\smallskip

Although the next result is essentially contained in
the paper \cite{WXL}, for completeness it is
worth stating and proving existence of (nonempty)
$\omega$-limit sets of weak solutions.
\bete\label{teo:omega}
 Let the assumptions of\/ {\rm Theorem~\ref{teo:esi}} hold,
 and let $(\bu,\bd)$ be a weak solution of {\rm Problem (P)}. Then, the $\omega$-limit set of $(\bu,\bd)$
 is nonempty. More precisely, we have
 \begin{equation}\label{conv:u}
   \lim_{t\nearrow+\infty} \bu(t) = 0
    \quext{weakly in }\,\bH,
 \end{equation}
 and any diverging sequence $\{t_n\}\subset[0,+\infty)$
 admits a subsequence, not relabeled, such that
 \begin{equation}\label{conv:d}
   \lim_{n\nearrow+\infty} \bd(t_n) = \bd\ii
    \quext{weakly in }\,\bV \quext{and strongly in }\,\bH,
 \end{equation}
 for some $\bd\ii\in V$. Moreover, any such limit point $\bd\ii$
 is a solution of the stationary problem
 \begin{equation}\label{eq:staz}
   - L \Delta \bz + \bbf(\bz) = 0~~\text{in }\,\Omega,
    \qquad \dn\bz=0~~\text{on }\,\Gamma.
 \end{equation}
\ente
\begin{proof}
 Let $\{t_n\}\subset[0,+\infty)$ be a diverging sequence.
 Then, the energy estimate implies that,
 at least for a (nonrelabeled) subsequence of $n$,
 \begin{equation}\label{omega:11}
   \bu(t_n)\to\bu\ii \quext{weakly in }\,\bH, \qquad
   \di(t_n)\to\di\ii \quext{weakly in }\,\bV,
 \end{equation}
 for suitable limit functions $\bu\ii$ and $\di\ii$.
 Let us consider the initial and boundary value problem
 associated to \eqref{eq:u}-\eqref{eq:d} on the time interval $[t_n,t_n +1]$ with ``initial
 values'' $\bu(t_n)$ and $\di(t_n)$. It is clear that,
 setting,  $\bu_n(t):=\bu(t+t_n)$ and
 $\di_n(t):=\di(t+t_n)$, $t\in[0,1]$, we get a weak solution
 to the problem on the time interval $[0,1]$. Then, \eqref{energy}
 implies that
 \begin{equation}\label{omega:12}
   \nabla\bu_n\to 0 \quext{strongly in }\,L^2(0,1;\bH^d),
 \end{equation}
 whence, by Poincar\'e's inequality and \eqref{energy} again,
 we have also
 \begin{equation}\label{omega:12b}
   \bu_n\to 0 \quext{strongly in }\,L^2(0,1;\bVz)
    \quext{and weakly star in }\,L^\infty(0,1;\bH).
 \end{equation}
 Moreover, we have
 \begin{equation}\label{omega:13}
   \di_n\to\ddi \quext{weakly star in }\,L^\infty(0,1;\bV)
    \cap L^2(0,1;H^2(\Omega)^d),
 \end{equation}
 for a suitable limit function $\ddi$.
 The growth condition \eqref{hp:psibound} and
 a comparison argument in \eqref{eq:d} then entail
 \begin{equation}\label{omega:14}
   \di_{n,t}\to\ddi_t \quext{weakly in }\,L^2(0,1;L^{3/2}(\Omega)^d).
 \end{equation}
 Hence, by the Aubin-Lions lemma, we obtain
 \begin{equation}\label{omega:14b}
   \di_{n}\to\ddi \quext{strongly in }\,L^2(0,1;\bV).
 \end{equation}
 To proceed, we take $\bphi \in \boldsymbol{W}^{1,3}_{0,\dive}(\Omega)$
 and test \eqref{eq:u}
 by $\bphi$. Noting that, by \eqref{eq:S},
 \begin{equation}\label{omega:15}
   \io\SSS_n:\nabla\bphi
    = - L \io ( \nabla \di_n \odot \nabla \di_n ) : \nabla \bphi
     - \delta \io \Big( \big( L \Delta \di_n - \bbf(\di_n) \big) \otimes \di_n \Big) : \nabla \bphi,
 \end{equation}
 and recalling \eqref{omega:12b}, \eqref{omega:13},
 and Remark~\ref{rem:weaksol}, we arrive at
 \begin{equation}\label{omega:17}
   \| \bu_{n,t} \|_{L^2(0,1;\boldsymbol{W}^{-1,3/2}_{\dive}(\Omega))}\le c,
 \end{equation}
 where $\boldsymbol{W}^{-1,3/2}_{\dive}(\Omega)$ denotes the dual space
 to $\boldsymbol{W}^{1,3}_{0,\dive}(\Omega)$ and
 $c$ denotes a positive constant independent of $n$.

 \smallskip

 Thus, from \eqref{omega:12b}, \eqref{omega:17}, and the Aubin-Lions lemma,
 we obtain  that
  \begin{equation}\label{hana222}
 \bu_n\to 0 \quad \mbox{strongly in } C^0([0,1];\bVz'),
 \end{equation}
 so that, in particular, $\bu\ii=0$, and \eqref{omega:12b}, \eqref{hana222} imply \eqref{conv:u}.
 On the other hand,
 by the energy estimate, we obtain
 \begin{equation}\label{omega:18}
   -L \Delta \di_n + \bbf(\di_n) \to 0
    \quext{strongly in }\,L^2(0,1;\bH),
 \end{equation}
 whereas, by \eqref{omega:12b}-\eqref{omega:14b},
 \begin{equation}\label{omega:19}
   \bu_n\cdot\nabla\di_n
    - \delta \di_n\cdot \nabla\bu_n \to 0
    \quext{weakly in }\,L^2(0,1;L^{3/2}(\Omega)^d).
 \end{equation}
 Thus, comparing terms in \eqref{eq:d},
 we also have that $\di_{n,t}\to 0$ in a suitable way.
 This entails that $\ddi$ is constant in time
 and, therefore, it coincides with $\di\ii$ for
 all times in $[0,1]$. Moreover,
 taking the limit in \eqref{eq:d}, we obtain that
 $\di\ii$ is a solution to \eqref{eq:staz}, as desired.
 This completes the proof.
\end{proof}
\noindent%
We now present the main results of this paper,
which characterize the $\omega$-limit set
of our system as a singleton under a number of different
conditions.
\bete\label{teo:1}
 Let the assumptions of\/ {\rm Theorem~\ref{teo:esi}} hold, and,
 in addition, let \/ $\psi$ be {\rm analytic}. Then the $\omega$-limit set of the component $\di$ of any
 weak solution consists of a single point, and
 we have
 \begin{equation}\label{conv:d:all}
   \lim_{t\nearrow+\infty} \bd(t) = \bd\ii
    \quext{strongly in }\, \bH
 \end{equation}
 for the whole trajectory $\bd$, where $\bd\ii$ is
 a solution to \eqref{eq:staz}.
\ente
\beos
As usually, when applying the \L ojasiewicz inequality, we can also get the rate of convergence of the form
 \[
 \|\di(t)-\di_\infty\|_{\bH}\le C(1+t)^{-\frac{\theta}{1-2\theta}},
 \]
 where $\theta$ is the \L ojasiewicz exponent, and $C$ is a suitably chosen
 constant depending on the initial energy and on the limit function.
 \eddos
 \noindent%
 In particular situations, we can also prove a stronger convergence result:
 \bete\label{teo:1.1}
 Under the hypotheses of\/ {\rm Theorem~\ref{teo:1}}, let, in addition,
 \begin{equation}\label{delta}
    \delta =0.
 \end{equation}
 Then,
 \begin{equation}\label{hana25}
   \bu(t)\to 0 \ \mbox{ strongly in }\,\bV,
 \end{equation}
 \begin{equation}\label{conv:d:all1}
   \lim_{t\nearrow+\infty} \bd(t) = \bd\ii
    \quext{strongly in }\, H^2(\Omega)^d.
 \end{equation}
\ente
\beos\label{periodic}
The same result was proved in \cite{WXL} for periodic boundary conditions
for $\bu$ and $\di$, large viscosity coefficient, and smooth initial data.
Actually, it is easy to check that our argument holds true also in
the case of periodic B.C., when $\delta\geq0$. Hence,
the same result of \cite{WXL} holds without the requirement of
large viscosity, and for initial data as in \eqref{hp:init1}, \eqref{hp:init2}.
On the other hand, in the case $\delta>0$ with boundary conditions
\eqref{bc:u}-\eqref{bc:d}, it does not seem possible to repeat the strong
estimates required for the proof of \eqref{hana25}-\eqref{conv:d:all1}
(some additionally boundary terms appear, which is not clear how to control).
Hence, extending the statement of Theorem~\ref{teo:1.1} to this situation
remains an open question.
\eddos
As in \cite[Thm.~1.2]{WXL} the proofs rely on a suitable
version of the {\sl Simon-\L ojasiewicz}\/ inequality,
proved in \cite[Thm.~6]{GG}. For the reader's convenience,
we report here the statement of a particular
case of the (more general) result of \cite{GG},
in a form suitable for our application:
\bete\label{teo:loj}
 Let the energy functional $\cal E$ be given by \eqref{defiE} with $\hat\psi$ analytic. Let $\pp\in \bV$ be a critical point of
 $\calE$. Then there exist constants $\theta\in(0,1/2)$,
 $\Lambda>0$ and $\epsilon_1>0$ such that the inequality
 \begin{equation}\label{ineq:loj}
   | \calE(\bv) - \calE(\pp) |^{1-\theta}
    \le \Lambda \big\| - L \Delta \bv + \bbf(\bv) \big\|_{\bV'}
 \end{equation}
holds for any $\bv$ such that
\begin{equation}\label{hana23}
\| \bv-\pp \|_{\bV}< \epsilon_1\,.
\end{equation}
\ente
To apply the preceding Theorem in our situation, we have to show that the inequality \eqref{ineq:loj}
holds for $\bv=\di(t)$ in a small $\bH$-neighbourhood of $\di_\infty$:
\bele\label{loj}
Let the energy functional $\cal E$ be given by \eqref{defiE} with $\hat\psi$ analytic.
Let $\di_\infty\in \bV$ be a solution of \eqref{eq:staz}. Let $K, \, P>0$ be constants.
 Then there exist $\epsilon>0$ and $\Lambda>0$ such that \eqref{ineq:loj} holds for any
 $\bv$ such that
 \begin{equation}\label{new}
 \|\bv\|_{\bV}\le K,\quad \| \bv-\di_\infty \|_{\bH}\le \epsilon, \hbox{ and}\quad  |\calE(\bv)-\calE(\di_\infty)|\le P\,.
\end{equation}
\enle
\begin{proof}
We argue by contradiction. Assume that there is a sequence $\bv_n$ such that
\[
\|\bv_n\|_{\bV}\le K,\ \ \bv_n\to \di_\infty \mbox{ in } \bH, \ \ |\calE(\bv_n)-\calE(\di_\infty)|\le P
\]
and
\begin{equation}\label{ineqn}
|\calE(\bv_n)-\calE(\di_\infty)|^{1-\theta}\ge n\|-L\Delta \bv_n+\bbf(\bv_n)\|_{\bV'},\ \ n=1,2,3,...
\end{equation}
Then
\[
\bbf(\bv_n)\to \bbf(\di_\infty)\mbox{ in } \bV', \mbox{ and }  \Delta \bv_n\to \Delta\di_\infty \mbox{ in } \bV'.
\]
This implies that
\[
\nabla\bv_n\to \nabla\di_\infty \mbox{ in } \bH, \mbox{ and, consequently, } \bv_n\to \di_\infty \mbox{ in } \bV.
\]
Hence, at least for $n$ sufficiently large, \eqref{hana23}
holds for $\bv=\bv_n$, $\pp=\di_\infty$. Consequently,
also \eqref{ineq:loj} is valid. This
contradicts \eqref{ineqn}.
\end{proof}

In the case that $\bbf$ does not satisfy the analyticity condition, we
can show that the $\omega$-limit set is a singleton
only in particular situations. For this purpose, we first
state a simple property:
\bele\label{lemma:gm}
 Let\/ \eqref{hp:psi}-\eqref{hp:psibound} hold. Then, $\ddi$ is
 a global minimizer of $\calE$ if and only if $\ddi$ is a
 constant unit vector.
\enle
\begin{proof}
 Thanks to \eqref{hp:psi}-\eqref{hp:psibound} the function
 $r\mapsto \psihat(r)-r$ has a minimum at $r=1$; moreover,
 $\psihat(1)-1=0$. Thus, $\calE(\di)$ is always nonnegative
 and $\calE(\di)=0$ if and only if $\nabla\di=0$ a.e.~in $\Omega$
 and $|\di|=1$ a.e.~in $\Omega$, whence the claim follows
 immediately.
\end{proof}
\noindent%
Our next result is
of conditional type and states that, if the $\omega$-limit
set of $\di(t)$ contains at least one global minimizer $\ddi$
of the free energy, then it has to coincide
with the set $\{\ddi\}$.
This is a consequence of the facts that the set of global
minimizers of the free energy is a $(d-1)$-dimensional smooth
manifold and, on the other hand, the kernel
of the linearized operator
$\bz \mapsto -\Delta \bz + \de_{\bd}\bbf(\ddi) \bz$
is also a $(d-1)$-dimensional manifold. In other words,
 the so-called {\sl normal hyperbolicity}\/
condition  is satisfied at $\ddi$, which
implies convergence of the whole trajectory to $\ddi$.
\bete\label{teo:2}
 Let the assumptions of\/ {\rm Theorem~\ref{teo:esi}} hold
 and let us assume that there exist a constant
 unit vector $\ddi\in\SSS^{d-1}$ and a diverging sequence
 $\{t_n\}$ such that
 \begin{equation}\label{conv:d:n}
   \lim_{t_n\nearrow+\infty} \bd(t_n) = \ddi
    \quext{weakly in }\,\bV.
 \end{equation}
 Then, $\omega-\lim\di=\{\ddi\}$
 and the whole trajectory $\di(t)$ converges to $\ddi$
 strongly in $\bH$ as $t\nearrow\infty$. If, in addition, \eqref{delta}  holds, then $\di(t)\to \ddi$ in $H^2(\Omega)^d$.
\ente

\beos\label{periodic2}
Let us note that the same convergence result for $\di$ in $H^2(\Omega)^d$ holds
true in case $\delta>0$ with periodic boundary conditions
for $\bu$ and $\di$.
\eddos

\noindent%
The next results only hold in the case
$\delta=0$. Actually, their proofs rely on the maximum
principle proved in Theorem~\ref{teo:max}.
In this setting, convergence to a single equilibrium takes place if either the
diffusion coefficient $L$ in \eqref{eq:d} is large
enough, or the ``initial energy'' $\EE_0:=\EE(\bu_0,\di_0)$
(cf.~\eqref{hp:init1}-\eqref{hp:init2}) is small enough (in
other words, if the initial datum $\di_0$ is
sufficiently close to the set of global minimizers). Indeed, we can prove the following two results:
\bete\label{teo:3}
 Let the assumptions of\/ {\rm Theorem~\ref{teo:max}}
 hold and, in particular,  let $\delta=0$. Assume that  $L$ in \eqref{eq:d} satisfies
 \begin{equation}\label{L}
  L>c_\Omega^2, \mbox{ where } c_\Omega
   \mbox{ is the best constant in the Poincar\'e-Wirtinger inequality}.\ \ \ \
 \end{equation}
 Then, the $\omega$-limit set of any weak solution
 starting from $(\bu_0,\di_0)$ consists
 of a single point $(\boldsymbol{0},\di_\infty)$.
\ente

\bete\label{teo:3bis}
 Let the assumptions of\/ {\rm Theorem~\ref{teo:max}}
 hold and, in particular,  let $\delta=0$. Assume that there exist $\kappa>0$ and $ \sigma\ge 1$ such that
 \begin{equation}\label{psi:1}
    \psihat(r)-r \ge \kappa (1-r)^\sigma
    \quad\perogni r\in[0,1].
 \end{equation}
 Then, there exists $\epsilon>0$ such that, if $(\bu_0,\di_0)$ satisfy $\EE_0 \le \epsilon$,
 then, the $\omega$-limit set of any weak solution
 starting from  $(\bu_0,\di_0)$ consists
 of a single point $(\boldsymbol{0},\di_\infty)$.
\ente


\section{Proofs}
\label{sec:proofs}

All proofs will be presented in the case $d=3$, the case $d=2$ being clearly
simpler.


%
%


\subsection{Proof of Theorem~\ref{teo:1}}
\label{subsec:1}

{\bf Energy estimate.}~~%
We test \eqref{eq:u} by $\bu$ and \eqref{eq:d} by
$-L \Delta \di + \bbf(\bd)$. Performing standard
computations and using, in particular, the incompressibility
constraint \eqref{incompr}, we readily obtain
the energy inequality~\eqref{energy}.
In particular,
we get that the function
$t\mapsto \EE(\bu(t),\di(t))$
is nonincreasing, whence it tends to
some (finite) value $\EE\ii$. Moreover,
thanks to \eqref{conv:u}-\eqref{conv:d},
we get
\begin{equation}\label{energy2}
  \EE\ii - \EE_0
   = - \int_0^{+\infty} \calD(s)\,\dis
   \le 0,
\end{equation}
where $\calD$ denotes the sum of the
{\sl dissipative terms}, namely
\begin{equation}\label{defiD}
  \calD := \big\| - L \Delta \di + \bbf(\di) \big\|^2
   + \nu \| \nabla \bu \|^2.
\end{equation}
We deduce from the energy inequality \eqref{energy} and \eqref{hp:psi}-\eqref{hp:psibound} that
 \begin{align}\label{reg:u1}
    \bu &\in L^\infty(0,\infty;\bH) \cap L^2(0,\infty;\bVz),\\
  \label{reg:d1}
   \di&\in  L^\infty(0,\infty;\bV), \\
  \label{reg:d3}
    -L\Delta \di+\bbf(\di) &\in L^2(0,\infty;\bH).
 \end{align}
Relations \eqref{reg:u1}, \eqref{reg:d1} imply
\begin{equation}\label{du}
\bu \cdot \nabla \di
   - \delta \di \cdot \nabla \bu\in L^2(0,\infty;L^{3/2}(\Omega)^d)
 \end{equation}
which, together with \eqref{reg:d3}, yield (cf. \eqref{eq:d})
\begin{equation}\label{dt}
\di_t\in L^2(0,\infty;L^{3/2}(\Omega)^d).
\end{equation}
\noindent%
{\bf Application of the \L ojasiewicz inequality}. Our aim is to show that there exists $T>0$ such that
\begin{equation}\label{dt1}
\di_t\in L^1(T,\infty;L^{3/2}(\Omega)^d),
\end{equation}
which implies convergence of $\di$ in $L^{3/2}(\Omega)^d$. The
pre-compactness of the trajectory
in $\bH$ then concludes the proof of Theorem \ref{teo:1}.

To this end, we first realize that there exists a constant $C$ such that
\begin{equation}\label{hana21}
\|\bu(t)\|^2\le C\|\nabla \bu(t)\|^{\frac1{1-\theta}} \ \mbox{ for a.a. } t>0,
\end{equation}
where $\theta\in(0,\frac12)$ is the same as in \eqref{ineq:loj}.
Indeed, if $\|\nabla \bu\|\le 1$, then \eqref{hana21}  follows by the Poincar\'e inequality,
and if  $\|\nabla \bu\|\ge 1$, the interpolation between $\bV_0$ and $\bV_0'$
together with the boundedness of $\bu$ in $\bV'$ gives the same estimate.

Now, let $\bd\ii$ be an element of the
$\omega$-limit set of $\bd$.
Then, integrating \eqref{energy} from 0 to
$+\infty$, we infer that
\begin{equation}\label{hana1}
\calD\in L^1(0,\infty),
\end{equation}
and, from Lemma~\ref{loj} and \eqref{hana21}, we get
\begin{equation}\label{hana11}
   \int_t^{+\infty}\calD(s)\,\dis = \calE(\bd(t)) - \calE(\bd\ii)
    + \frac12 \| \bu(t) \|^2 \le C \calD(t)^{\frac1{2(1-\theta)}},
\end{equation}
for all $t>0$ such that \eqref{new} holds. Denoting
by $\calM$ this set,
 we obtain that $\calD^{1/2} \in L^1(\calM)$, see \cite[Lemma~7.1]{FS}. Then also
\begin{equation}\label{hana12}
   \bu \in L^1(\calM;\bVz),~~
   - L \Delta \di + \bbf(\di) \in L^1(\calM;\bH),
\end{equation}
and so, taking into account the growth of $\bbf$, we get
\[
\bu\cdot\nabla\di
    - \delta \di\cdot \nabla\bu\in L^1(\calM,L^{\frac32}(\Omega)^d),
    \]
which implies (cf.~\eqref{eq:d})
\begin{equation}\label{hana22}
\di_t\in L^1(\calM,L^{\frac32}(\Omega)^d).
\end{equation}
This fact, combined with the pre-compactness of the trajectory of $\di$ in $\bH$ and a simple contradiction
argument (see \cite{FS}), yields the existence of
a large $T$ such that (cf. Lemma~\ref{loj})
\[
\|\di(t)-\di_\infty\|<\epsilon, \quad\forall t\geq T.\ \
\]
In other words, the solution $\di$
remains in the $\epsilon-$neighbourhood of $\di_\infty$ in the space $\bH$ for $t\ge T$,
and Lemma \ref{loj} applies to $\di(t)$ in the whole interval $(T,\infty)$. In other words, ${\cal M}\supset (T,\infty)$,
and \eqref{dt1} holds. This fact, together with pre-compactness of the trajectory yields
\begin{equation}\label{hana15}
  \di(t) \to \di\ii \quext{strongly in }\,\bH.
\end{equation}
Theorem \ref{teo:1} has been proved.

\subsection{Proof of Theorem~\ref{teo:1.1}}
\label{subsec:1.1}

The proof follows the lines of the argument developed
in \cite[Sec.~3.2]{CGR}. In particular, we can
derive the differential inequality for the dissipative
term $\calD$ defined in \eqref{defiD} (cf. \cite[formula (8)]{CGR}, which is still valid with our
boundary conditions in case $\delta=0$ or with periodic boundary conditions in case $\delta>0$):
\begin{equation}\label{ineq:WXL}
  \ddt \calD \le C_* \big( \calD^3 + 1 \big),
\end{equation}
where the computable constant $C_*$ depends
on the parameters of the problem and on the ``initial energy''
$\EE_0$, but is independent of time.
We point out that the formal computations in the proof
of \eqref{ineq:WXL} are valid for the classical solutions, but the result can be justified by a proper approximation.

 Next, we show that there is $T_1>0$ such that $\calD\in L^\infty(T_1,\infty)$.
 To prove this,
we consider the differential inequality
\begin{equation}\label{ineq:WXL-2}
  y' \le C_* \big( y^3 + 1 \big),
   \qquad y(t_0)=1.
\end{equation}
Then, there exist (a small) $\tau$
(independent of $t_0$)
and (a large) $K>0$ such that the solution
$y$ satisfies
\begin{equation}\label{ineq:WXL-3}
  \| y \|_{C^0([t_0,t_0+\tau])} \le K.
\end{equation}
On the other hand, according to \eqref{hana1},
we have
\begin{equation}\label{conto0-11}
  \lim_{t\nearrow+\infty} \int_t^{+\infty}
   \calD(s)\,\dis = 0.
\end{equation}
Thus, for any $\epsilon>0$ there exists $T>0$ such that
\begin{equation}\label{conto0-12}
  \int_T^{+\infty} \calD(s)\,\dis \le \epsilon.
\end{equation}
Choosing $\epsilon=\tau/2$ and $T$ correspondingly,
we obtain that, for all $t\ge T$,
there exists $t_0\in [t,t+\tau/2]$
such that
\begin{equation}\label{conto0-13}
  \calD(t_0)
   \le \frac{2}\tau \int_t^{t+\tau/2} \calD(s)\,\dis
   \le \frac{2\epsilon}\tau
  = 1.
\end{equation}
Comparing solutions of
\eqref{ineq:WXL} and \eqref{ineq:WXL-2} and recalling the choice
of $t_0$, we get from \eqref{ineq:WXL-3} that
\begin{equation}\label{ineq:WXL-3b}
  | \calD(s) | \le K
   \quad \perogni s\in [T+\tau/2,+\infty).
\end{equation}
Setting $T_1:=T+\tau/2$, we deduce from \eqref{ineq:WXL} that
$\ddt\calD$ is bounded on $(T_1,\infty)$, which together with \eqref{hana1} yields
\begin{equation}\label{hana2}
\calD(t)\to 0 \ \mbox{ as } t\to \infty.
\end{equation}
This implies (using the Poincar\'e inequality) that,
\begin{equation}\label{couv}
   \bu(t) \to 0, \quext{strongly in }\,\bV.
\end{equation}
Taking into account the growth of $\bbf$, we get from \eqref{ineq:WXL-3b}
\begin{equation}\label{regoult:d}
 \| \bd \|_{L^\infty(T_1,\infty;H^2(\Omega)^d)} \le K.
\end{equation}
To show that $\bd$ converges to a single point $\bd_\infty$,
we make again use of the \L ojasiewicz inequality \eqref{ineq:loj}.
The same argument applies this time to the strong solution and time $t\ge T_1$.
This gives
\begin{equation}\label{hana14}
  \di_t \in L^1(T,\infty;\bH)\ \mbox{ for some }T>T_1.
\end{equation}
It follows that
\begin{equation}\label{hana15bis}
  \di(t) \to \di\ii \quext{strongly in }\,\bH.
\end{equation}
Moreover, by  \eqref{regoult:d} and the growth conditions on $\bbf$,
\begin{equation}\label{hana16}
  \bbf(\bd(t))\to\bbf(\bd_\infty)\  \mbox{ strongly in } \bH,
\end{equation}
whence, using \eqref{hana2} and \eqref{couv},
\begin{equation}\label{hana17}
  \Delta\bd(t)\to\Delta \bd_\infty\ \mbox{ strongly in } \bH,
\end{equation}
which concludes the proof.


\subsection{Proof of Theorem~\ref{teo:2}}
\label{subsec:2}

In this section, we show that the energy functional $\cal E$ satisfies
the \L ojasiewicz inequality \eqref{ineq:loj} with the exponent $\frac12$.
Then, arguing as in the proof of Theorem \ref{teo:1}, we obtain the strong convergence
of $\di$ in $\bH$. Moreover, if \eqref{delta} holds (or we have periodic boundary
conditions, cf.~Remark~\ref{periodic}),
we have
the strong convergence in $H^2(\Omega)^d$ (cf. the proof of Theorem~\ref{teo:1.1}).

Let us consider the linearized problem associated to
\eqref{eq:staz} at the element $\ddi$ of the
$\omega$-limit set, i.e.,
\begin{equation}\label{eq:staz-lin}
  \calL(\ddi)\bz:= - L \Delta \bz
  + \psi(|\ddi|^2)\bz
  + 2 \psi'(|\ddi|^2)(\ddi\otimes\ddi) \bz
  - \bz = 0,
   \qquad \dn\bz=0~~\text{on }\,\Gamma.
\end{equation}
Let $\ddi$ be a global minimizer
of $\calE$, i.e., a constant unit vector
(by Lemma~\ref{lemma:gm}).
We aim to apply the result proved by Simon and reported in \cite[Cor. 3.12]{Chi}.
To this end, we introduce the following notation:
\[
U \mbox{ is a $\bV$-neighbourhood of } \ddi \in \bV,
\]
\[
{\boldsymbol{\cal V}}_0 \mbox{ is the kernel of } \ \calL(\ddi),
\]
\[
S_0=\{ \bd\in \bV ;\  \calE'(\bd)=0\},
\]
\[
S=\{ \bh\in U;\ \calE'(\ddi+\bh)\in {\boldsymbol{\cal V}}_0'\}.
\]
In our situation, \cite[Cor. 3.12]{Chi} reads as follows:
\bele\label{Simon}
Let $\ddi\in S_0$ and assume the following hypotheses:

(i) The kernel ${\boldsymbol{\cal V}}_0$ of the linearization $ \calL(\ddi)$ is
a complemented subspace of  $\bV$, i.e., there exists a projection $P\in {\calB}(\bV)$ such that ${\boldsymbol{\calV}}_0={\rm Rg} P$.

(ii) There exists a neighbourhood $U$ of $\ddi$ in $\bV$ such that $\calE'\in C^1(U,\bV')$.
Moreover, the range of $\calL(\ddi)$ coincides with  $\boldsymbol{\cal V}_1'$, the space of the
elements of $\boldsymbol{V}'$ belonging to the the kernel of the adjoint projection $P'\in {\cal B}(\bV')$.

(iii) $(S_0-\ddi)\cap S$ is a neighbourhood of 0 in the critical manifold $S$.

\noindent
Then $\cal E$ satisfies the \L ojasiewicz inequality near $\ddi$ with the exponent $\theta=\frac12$.
\enle
\noindent%
To verify the assumption (i), we  test \eqref{eq:staz-lin} by
$\bz$ and use the condition $\psi(1)=1$ to obtain
\begin{equation}\label{co:2-01}
  L \| \nabla \bz \|^2
  + 2 \io \psi'(1) | \ddi \cdot \bz |^2
  = 0.
\end{equation}
Hence, taking into account  the last condition in \eqref{hp:psi}, we get
\begin{equation}\label{co:2-02}
  \nabla \bz = 0 \quand
   \ddi \cdot \bz = 0
   \quext{a.e.~in }\,\Omega.
\end{equation}
Consequently, any solution $\bz$ to \eqref{eq:staz-lin},
i.e., any element of the kernel, is
a constant vector orthogonal to $\ddi$ (conversely,
it is apparent that any such vector is a solution to
\eqref{eq:staz-lin}). Thus, the kernel of the
linearized operator $\calL(\ddi)$
is a $(d-1)$-dimensional  plane orthogonal
to $\ddi$ and containing the origin, which trivially
permits to define the projection $P$.\\[2mm]
The first condition in (ii) is obvious since $\bbf$ is
$C^1$ and, by hypotheses, has at most cubic growth.
%
%
To verify the second condition, we observe that $\boldsymbol{\calV}_1'$ is the
subspace of $\bV'$ consisting of the elements that are orthogonal
(w.r.t.~the duality between $\bV'$ and $\bV$) to the plane $\boldsymbol{\calV}_0$.
Then, computing $\duav{\calL(\ddi)\bz,\bv_0}$ for generic $\bz\in \bV$
and $\bv_0\in \boldsymbol{\calV}_0$, we obtain (cf.~\eqref{eq:staz-lin} and recall
that $\psi(1)=1$), using \eqref{co:2-02},
\begin{equation}\label{orthog}
  \duav{\calL(\ddi)\bz,\bv_0}
  = 2 \io \psi'(|\ddi|^2) (\ddi \cdot \bz)
   (\ddi \cdot \bv_0) = 0,
\end{equation}
the last equality following from the fact that
$\bv_0\perp\ddi$. Thus, $\calL(\ddi)\bV\subset \boldsymbol{\calV}_1'$.
To show the converse inclusion, we choose $\boldsymbol{\zeta}\in \boldsymbol{\calV}_1'$
and prove that there exists at least one $\bz \in \bV$ such that
\begin{equation}\label{surj}
  \calL(\ddi)\bz = - L \Delta \bz
   + 2 \psi'(1)(\ddi\otimes\ddi) \bz
   = \boldsymbol{\zeta},
   \qquad \dn\bz=0~~\text{on }\,\Gamma.
\end{equation}
This can be seen by approximation. Actually, it is clear that,
for any $k\in\NN$, there is a solution $\bz_k$ to
\begin{equation}\label{surj2}
  - L \Delta \bz_k
  + k^{-1} \bz_k
   + 2 \psi'(1)(\ddi\otimes\ddi) \bz_k
  = \boldsymbol{\zeta},
   \qquad \dn\bz=0~~\text{on }\,\Gamma.
\end{equation}
Testing by $\bz_k$, we have
\begin{equation}\label{surj3}
  L \| \nabla \bz_k \|^2
  + k^{-1} \| \bz_k \|^2
   + 2 \psi'(1) \| \ddi \cdot \bz_k \|^2
  = (\boldsymbol{\zeta}, \bz_k) = (\boldsymbol{\zeta}, \bz_k - P\bz_k).
\end{equation}
Indeed, $\boldsymbol{\zeta}\perp P\bz_k$ by assumption. Using the
Poincar\'e-Wirtinger inequality it is then apparent
that the \rhs\ can be estimated. Then, standard methods
permit to check that $\bz_k$ tend to a solution $\bz$
to \eqref{surj}, as desired.\\[2mm]
The third assumption is satisfied because
$
0\in S_0-\ddi\subset S,
$
and both $ S_0-\ddi$ and $S$ have the same dimension (see \cite[Proposition 3.6]{Chi}).

Lemma \ref{Simon} then yields that the
\L ojasiewicz inequality \eqref{ineq:loj}
holds near $\ddi$ with the exponent $\theta=1/2$.
If $\delta=0$ or in case of periodic boundary conditions, repeating the computations leading to \eqref{regoult:d},
we obtain the strong convergence of $\di$ in $\bH$, and,
proceeding as in \eqref{hana16}, \eqref{hana17},
the strong convergence in $H^2(\Omega)^d$,
which completes the proof of Theorem \ref{teo:2}.


\subsection{Proof of Theorem~\ref{teo:3}}
\label{subsec:3}

Let $\dii$ be an element of the $\omega$-limit set
of the $\di$-component of some weak solution. Then,
by Theorem~\ref{teo:omega},
$\dii$ solves the stationary problem, which
we rewrite as
\begin{equation}\label{prostaz1}
  - L \Delta \dii + \big( \psi (|\dii|^2) - 1 \big) \dii = 0.
\end{equation}
 Testing \eqref{prostaz1} by $\dii-(\dii)\OO$, where $(\dii)\OO=\int_\Omega\dii$, we get
\begin{equation}\label{contostaz11}
  L \| \nabla \dii \|^2 + \io \psi (|\dii|^2) \dii \cdot \big(\dii - (\dii)\OO \big)
   = \big\| \dii - (\dii)\OO \big\|^2 \le c\OO^2 \| \nabla \dii \|^2,
\end{equation}
where $c\OO$ is the (best) constant in the Poincar\'e-Wirtinger
inequality. Thus, being  $L$ is large (precisely,
we need $L>c\OO^2$), the latter term can be controlled.

In what follows, we denote
\begin{equation}\label{prostaz11}
  \Xi(\di)=\frac12 \psihat(|\di|^2).
\end{equation}
A direct check (e.g., computing the Hessian matrix)
shows that $\Xi$ is convex, and
we notice that $\de_{\di}\Xi(\di)=\psi (|\di|^2) \di$.
Thus, we have
\begin{align}\no
  \io \psi (|\dii|^2)\dii\cdot \big(\dii - (\dii)\OO \big)
   & = \big( \de_{\di}\Xi(\dii) , \dii - (\dii)\OO \big)_H\\
 \no
   & \ge \io \big( \Xi(\dii) - \Xi((\dii)\OO) \big)\\
 \label{contostaz12}
  & = \io \Xi(\dii) - \Xi\Big(\io \dii\Big)
    \ge 0,
\end{align}
where the latter inequality follows from Jensen's inequality.
%
%
%
%
%
%
%
%

\smallskip

From \eqref{contostaz11}-\eqref{contostaz12},
we obtain that $\dii$
is a constant vector. Thus, taking into account that $\psi$ is monotone
and $\psi(1)=1$, we readily obtain from equation \eqref{prostaz1}
 that either $\dii=0$ or $|\dii|=1$. Hence the set of stationary solutions
is disconnected and consists of the isolated point
and the two-dimensional manifold. Consequently,
either Theorem~\ref{teo:2} applies, or the whole
trajectory tends to $0$. In both
cases, the $\omega$-limit set is a singleton, as
desired.


\subsection{Proof of Theorem~\ref{teo:3bis}}
\label{subsec:3bis}

Let us first note that,
by \eqref{energy2} and \eqref{psi:1},
we have
\begin{equation}\label{conto3-11}
  \epsilon \ge
  \EE_0 \ge \EE\ii
  = \frac12 \io \big(L |\nabla \dii|^2
    + \psihat(|\dii|^2)
    - |\dii|^2 \big).
\end{equation}
Rewriting the stationary problem  \eqref{prostaz1}
and testing it by $\dii$, we obtain
\begin{equation}\label{conto3-12}
  L \|\nabla \dii\|^2
   + \io \big(\psi(|\dii|^2)|\dii|^2-|\dii|^2\big)
   = 0.
\end{equation}
Dividing \eqref{conto3-12} by $2$ and subtracting
the result from \eqref{conto3-11}, we obtain
\begin{equation}\label{conto3-13}
  \frac12\io \big( \psihat(|\dii|^2)
   - \psi(|\dii|^2)|\dii|^2 \big)
   \le \epsilon.
\end{equation}
On the other hand, thanks to \eqref{hp:psihat}, \eqref{hp:psi}
and \eqref{psi:1},
\begin{equation}\label{co3-14}
  \frac12 \io \big( \psihat(|\dii|^2)
   - \psi(|\dii|^2)|\dii|^2 \big)
  \ge  \frac12 \io \big( \psihat(|\dii|^2)
   - |\dii|^2 \big)
  \ge \frac\kappa2 \io \big| 1 -|\dii|^2 \big|^\sigma,
\end{equation}
where also the maximum principle \eqref{maxprinc}
has been used. Thus,
\begin{equation}\label{conto3-15}
  \big\| 1 - |\dii|^2 \big\|_{L^1(\Omega)}
   \le \big\| 1 - |\dii|^2 \big\|_{L^\sigma(\Omega)}
   \le \Big(\frac{2\epsilon}{\kappa}\Big)^{1/\sigma}.
\end{equation}
To proceed, we notice that, by standard elliptic
regularity results applied to \eqref{eq:staz},
there exists a constant $K_0>0$ such that,
for any solution $\ddi$ of \eqref{eq:staz}
it holds
\begin{equation}\label{conto3-16}
  \big\| \ddi \big\|_{H^2(\Omega)} \le K_0.
\end{equation}
Consequently, for some $K>0$ depending on $K_0$,
we have
\begin{equation}\label{conto3-17}
  \big\| 1 - |\ddi|^2 \big\|_{H^2(\Omega)} \le K.
\end{equation}
In particular, $\di\ii$ satisfies \eqref{conto3-17}.
Thus, by the Gagliardo-Nirenberg interpolation
inequality (we refer, for simplicity, to the case
$d=3$, the case $d=2$ is even better),
\begin{align}\no
  \big\| 1 - |\dii|^2 \big\|_{L^\infty(\Omega)}
   & \le \big\| 1 - |\dii|^2 \big\|_{H^2(\Omega)}^{6/7}
      \big\| 1 - |\dii|^2 \big\|_{L^1(\Omega)}^{1/7}
     + \big\| 1 - |\dii|^2 \big\|_{L^1(\Omega)}\\
 \label{conto3-18}
  & \le K^{6/7}\Big(\frac{2\epsilon}{\kappa}\Big)^{1/7\sigma}
   + \Big(\frac{2\epsilon}{\kappa}\Big)^{1/\sigma}
  < \eta,
\end{align}
where $\eta>0$ is a small constant to be chosen later
and the last inequality holds provided that
$\epsilon$ is small enough.

To conclude the proof, we set
\begin{equation}\label{defialfa}
  \alpha:=\psi(|\dii|^2)-1
\end{equation}
and notice that $\alpha\leq 0$ because of \eqref{hp:psi}-\eqref{hp:psihat}. Moreover,
$\dii$ can be interpreted as a solution
of the linear elliptic system
\begin{equation}\label{ell:lin}
  - L \Delta \dii + \alpha \dii = 0~~\text{in }\,\Omega,
   \qquad \dn\dii=0~~\text{on }\,\Gamma.
\end{equation}
Testing the above equation by $1$, we obtain
\begin{equation}\label{conto3-20}
  \io \alpha \dii = 0.
\end{equation}
Then, multiplying \eqref{ell:lin} by $\dii-(\dii)\OO$
and using \eqref{conto3-20}, we infer
\begin{align}\no
  L \| \nabla \dii \|^2
   & = -\io \alpha \big| \dii - (\dii)\OO \big|^2
    - (\dii)\OO \cdot \io \alpha \big( \dii - (\dii)\OO \big)\\
 \label{conto3-21}
   & \le c\OO^2 \| \alpha \|_{L^\infty(\Omega)}\| \nabla\dii \|^2
    + \big| (\dii)\OO \big|^2 \io \alpha.
\end{align}
The latter term is nonpositive, while the first
term on the \rhs\ can be controlled provided that $\eta$ is small
enough. Indeed,
\begin{align}\no
  \| \alpha\|_{L^\infty(\Omega)}
   & = \big\| \psi(|\dii|^2)-1 \big\|_{L^\infty(\Omega)}
   = \big\| \psi(|\dii|^2)-\psi(1) \big\|_{L^\infty(\Omega)}\\
 \label{conto3-22}
   & \le c_\psi \big\| |\dii|^2 - 1 \big\|_{L^\infty(\Omega)}
    \le c_\psi \eta,
\end{align}
thanks to \eqref{conto3-18} and \eqref{hp:psibound}. If $L>c_\Omega^2 c_\psi\eta$, we see, in the same way as above, that
$\dii$ is a constant unit vector. Finally, we take $\epsilon$ such that \eqref{conto3-18} holds, and
the proof follows
again by applying Theorem~\ref{teo:2}.




\end{document}